\documentclass[12pt]{article}
\usepackage{amssymb,amsmath,amsthm,fancyhdr,array,calc}
\let\x=\times
\def\HW{{\textsc{\small HW}\!}}
\def\ie{{\it i.e.}}

\def\OpGets#1{{\;#1\!\!:=\;}}
\normalsize
\headwidth=\textwidth
\newcommand{\fancytitle}{\maketitle\thispagestyle{fancy}}

\title{\small\bfseries \Large\bfseries Gray-coding through nested sets }
\author{Antonia W.~Bluher, National Security Agency }
\date{February 9, 2015}

\begin{document}
\fancytitle
\begin{abstract}
 We consider the following combinatorial question. Let
$$ S_0 \subset S_1 \subset S_2 \subset \ldots \subset S_m $$
be nested sets, where $\#(S_i) = i$. A {\it move} consists of altering one 
of the sets $S_i$, $1 \le i \le m-1$, in a manner so that the nested 
condition still holds and $\#(S_i)$ is still $i$. Our goal is to find a 
sequence of moves that exhausts through all subsets of $S_m$ (other than 
the initial sets $S_i$) with no repeats. We call this ``Gray-coding through 
nested sets'' because of the analogy with Frank Gray's theory of 
exhausting through integers while altering only one bit at a time. 
Our main result is an efficient algorithm that solves this problem.
As a byproduct, we produce new families of cyclic Gray codes through
binary $m$-bit integers.
\end{abstract}

\section{Introduction}

The nested set problem can be stated as follows.  Let
$$ S_0 \subset S_1 \subset S_2 \subset \ldots \subset S_m $$
be nested sets, where $\#(S_i) = i$. A {\it move} consists of altering one 
of the sets $S_i$,
for $1\le i < m$, while maintaining the conditions that $\#(S_i)=i$ and that
the sets are nested. Our goal is to find a sequence of moves that exhausts
through all subsets of $S_m$ (other than the initial sets $S_i$) with
no repeats. We call this ``Gray-coding
through nested sets'' because of the analogy with Frank Gray's theory
on stepping through all $m$-bit integers by altering just one bit at a time.
(For the theory of Gray codes see \cite{Gray} or
\cite[Section 2.2.2]{KreherStinson}.)
The problem of Gray-coding through nested sets does not seem
to have been considered before in the literature, but it fits naturally
into the genre of ``combinatorial Gray codes'' that has been 
popularized by Hilbert S.~Wilf, Carla Savage, F.~Ruskey, and others.
We quote from Savage's survey article \cite{Savage}: 
``The term combinatorial Gray code
$\dots$ is now used to refer to any method for generating combinatorial objects 
so that successive objects differ in some prespecified, usually small, way.''
Examples of such combinatorial objects include $k$-element subsets of
an $n$-element set \cite{BW,Chase,Wilf}, permutations, binary trees, or partitions.

Given nested sets $S_i$ as above, define the difference sequence to be 
$\langle q_1,q_2,\ldots,q_m \rangle$,
where $q_i$ is the unique element of $S_i \setminus S_{i-1}$. 
The difference sequence gives enough information to determine all
the sets $S_i$.
The only valid way to alter $S_i$ is to remove $q_i$
and replace it with $q_{i+1}$, and this has the effect of transposing
$q_i$ and $q_{i+1}$ in the difference sequence. 
Then $S_i$ is replaced by the set $\{q_1,q_2,\ldots,q_{i-1},q_{i+1}\}$.

Since a move is completely determined by~$i$, we can record
 a sequence of moves by listing the indices.
That is, the sequence $[i_1,i_2,i_3,\ldots,i_N]$ means
to first alter $S_{i_1}$, then $S_{i_2}$, {\it etc.}
We call this a \emph{stepping sequence}
if every subset of $S_m$ (other than the initial sets) appears exactly once.
In other words, a stepping sequence is a sequence of moves that solves
the problem of Gray-coding through nested sets.
Since there are $2^m$ subsets and $m+1$ initial sets $S_i$,
a stepping sequence for $m$ has length $2^m-m-1$.
If $S_i$ is altered twice in a row, the set returns to its
original value.  This is not allowed since we are not supposed to
generate a set that has been seen before. Thus, any two consecutive terms
in a stepping sequence must be distinct. Also notice that $i$ occurs
in the stepping sequence exactly $\left({m \atop i}\right)-1$ times.

We first consider small examples. For $m=1$, we have $S_0=\emptyset \subset
S_1 = \{q_1\}$.  There are no additional subsets of $\{q_1\}$ other than
the initial sets $S_0$ and $S_1$, so the stepping sequence is empty.
For $m=2$, we begin with
$$S_0 = \emptyset \subset S_1 = \{q_1\} \subset S_2 = \{q_1,q_2\}.$$
The only subset of $S_2$ that has not been seen is $\{q_2\}$. 
We need a single move, which is to alter $S_1$,
and the associated stepping sequence is $[1]$.

For $m=3$ we begin with
$$S_0 = \emptyset \subset
S_1 = \{q_1\} \subset S_2 = \{q_1,q_2\} \subset S_3 = \{q_1,q_2,q_3\}.$$
One solution is as follows: change $S_2$ to $\{q_1,q_3\}$, 
change $S_1$ to $\{q_3\}$, change $S_2$ to $\{q_2,q_3\}$,
then change $S_1$ to $\{q_2\}$.  This is codified by $[2,1,2,1]$.
The only other solution is to first change $S_1$ to $\{q_2\}$,
then change $S_2$ to $\{q_2,q_3\}$, then change $S_1$ to $\{q_3\}$,
then change $S_2$ to $\{q_1,q_3\}$.   This solution is codified by
$[1,2,1,2]$.

Since the sets $S_i$ are nested, the complements $S_m \setminus S_i$
are also nested, but with the inclusions going in the reverse direction.
When we alter $S_i$, we also alter its complement.  Further, if the $S_i$'s
run through every possible subset, then so do the complements.
The complement of $S_i$
has order $m-i$.  Thus, if $[i_1,\ldots,i_N]$ is a stepping
sequence for $m$, then so is $[m-i_1,\ldots,m-i_N]$.  We call this the
{\it complementary stepping sequence}, or more simply, the complement.
In particular, for the order-3 example above,
the complement of $[2,1,2,1]$ is $[1,2,1,2]$.

In addition, the reverse of a stepping sequence is a stepping sequence.
This can be seen by running a stepping sequence backwards
(beginning with the final set values $S_i'$) and noting that
one finishes with the sets $S_i$, and in the intermediate steps one
passes through all remaining subsets of $S_m$.  The stepping
sequence $[2,1,2,1]$ has the property that its complement is equal to its
reverse, but not all stepping sequences have that property.
For $m=4$, there are exactly 34 stepping sequences, and ten of these
have the property that the complement is equal to the reverse.

In this article, we present four ways to create stepping sequences:
by recursion, by greed, with a for-loop, and with a different for-loop.
It turns out that all four methods are equivalent, {\it i.e.}, they give
rise to the same stepping sequence.  In the final section, we show
how the nested set problem is related to other combinatorial Gray codes,
and we pose some open problems.

\section{The recursive method}

The following theorem enables us to build stepping
sequences recursively.

\outer\def\proclaim#1. #2 \par{\medbreak\noindent{\bf #1.\enspace}{\it #2}\par
\ifdim\lastskip<\medskipamount\removelastskip\penalty55\medskip\fi}
\newenvironment{pf}{\noindent {\bf Proof.}}{\qed\vskip 6pt}

\proclaim{Theorem 2.1}. {If $[i_1,\ldots,i_K]$ and $[i_1',\ldots,i_K']$ are
stepping sequences for $m-1$ (so $K=2^{m-1}-m$), then
\begin{equation}
[i_1+1,i_2+1,\ldots,i_K+1] \cup [1,2,\ldots,m-1] \cup [i'_1,i'_2,\ldots,i'_K]
\label{recursion1}
\end{equation}
and
\begin{equation}
[i_1,i_2,\ldots,i_K] \cup [m-1,\ldots,1] \cup [i'_1+1,i'_2+1,\ldots,i'_K+1]
\label{recursion2}
\end{equation}
are stepping sequences for $m$. }

\begin{proof}  We first prove that (\ref{recursion1})
is a stepping sequence for $m$.
Suppose initially $S_1= \{a\}$.
Let $S_i' = S_i \setminus \{a\}$ for $i=1,\ldots,m-1$.
If we applied $[i_1,\ldots,i_K]$ to $S_i'$, we would obtain
every subset of $S_m'$, except for the initial sets $S_i'$.
By doing the same sequence of moves, but on the sets $S_i=S_i' \cup \{a\}$
instead of the sets $S_i'$,
we run through every subset of $S_m$ that contains $a$, except for
the initial sets $S_1,\ldots,S_m$.  Since $|S_i|=|S_i'|+1$,
this sequence of moves is codified as $[i_1+1,\ldots,i_K+1]$.
At this point we have seen exactly once every subset of $S_m$ that 
contains~$a$, and no other nonempty subset, and the difference sequence
has the form $\langle a,q_1,q_2,\ldots,q_{m-1} \rangle$.
Note that $\{q_1,q_2,\ldots,q_{m-1}\} = S_m \setminus \{a\}$.

Now we have nested sets
$$S_0 = \emptyset,\quad  S_1 = \{a\}, \quad S_2 = \{a,q_1\}, \quad
S_3 = \{a,q_1,q_2\}, \ \ldots, \  S_m=\{a,q_1,\ldots,q_{m-1}\}.$$
The moves $[1,2,\ldots,m-1]$ change $S_1$ to $\{q_1\}$, then change
$S_2$ to $\{q_1,q_2\}$, then change $S_3$ to $\{q_1,q_2,q_3\}$, and so on.
All these subsets are new, because they do not contain $a$. After these moves,
we have $S_i = \{q_1,q_2,\ldots,q_i\}$ for $i<m$. In particular, 
$S_{m-1}=\{q_1,\ldots,q_{m-1}\} = S_m \setminus \{a\}$.

Now applying the moves $[i'_1,i'_2,\ldots,i'_K]$ to the nested sets 
$S_0,S_1,\ldots,S_{m-1}$
produces all remaining subsets of $S_{m-1}$, {\it i.e.}, all remaining subsets
of $S_m$ that do not contain $a$.

In summary,
after applying the moves $[i_1+1,\ldots,i_K+1]$ we have seen all sets that contain $a$
exactly once.
The moves $[1,2,\ldots,m-1]$ produce
new sets $S_1=\{q_1\}$, $S_2=\{q_1,q_2\}$, $\ldots$, $S_{m-1}=\{q_1,\ldots,q_{m-1}\}$
that do not contain $a$.
Finally, the moves $[i_1,i_2,\ldots,i_K]$ produce
all remaining subsets that do not contain $a$, and each set is produced exactly once.
Thus all subsets of $S_m$ are produced exactly once, proving that 
(\ref{recursion1}) is indeed a stepping sequence for~$m$.

By applying (\ref{recursion1}) to the complements $[m-1-i_1,\ldots,m-1-i_K]$
and $[m-1-i_1',\ldots,m-1-i_K']$, we see that
$$[m-i_1,m-i_2,\ldots,m-i_K] \cup [1,2,\ldots,m-1] \cup [m-1-i_1',m-1-i_2',\ldots,m-1-i_K']$$
is a stepping sequence for $m$.  Taking the complement, we find that
(\ref{recursion2}) is also a stepping sequence for $m$.
\end{proof}

The stepping sequences for $m$ that are produced in the theorem
have length $2K+m-1=2(2^{m-1}-(m-1)-1) +m-1 = 2^m-m-1$, as expected.

The theorem enables one to build stepping
sequences for $m$ out of smaller stepping sequences. In particular, we define:
\begin{equation} \label{Rm}
\text{$R_2=[1]$, and $R_{m} = (R_{m-1} + 1) \cup [1,2,\ldots,m-1] \cup R_{m-1}$ for
$m=3,4,\ldots$,}
\end{equation}
where $R_{m-1}+1$ is obtained by adding one to each element of $R_{m-1}$.
Then $R_m$ is a stepping sequence for $m$ by Theorem~2.1.
The first few are given by:
\begin{eqnarray*}
R_2 &=& [1] \\
R_3 &=& [2]\cup[1,2]\cup [1] = [2,1,2,1] \\
R_4 &=& [3,2,3,2] \cup [1,2,3] \cup [2,1,2,1] = [3,2,3,2,1,2,3,2,1,2,1] \\
R_5 &=& [4,3,4,3,2,3,4,3,2,3,2,1,2,3,4,3,2,3,2,1,2,3,2,1,2,1].
\end{eqnarray*}

Figure~1 below illustrates nested-set generation for the stepping
sequence $R_4$ when $S_i=\{1,2,\ldots,i\}$. The stepping sequence is 
listed in the first
column, the difference sequence in the second column, and the
newly generated sets in the remaining columns.
Notice that the sets containing 1 are generated first, then the
sets not containing 1.  

\begin{center}
\begin{figure}[h]
\begin{tabular}{|c|c|c|c|c|c|}
\hline
Stepping & Difference & $S_1$ & $S_2$ & $S_3$ & $S_4$ \\
Sequence & Sequence &&&& \\
\hline
Initial & $\langle 1,2,3,4 \rangle$ & $\{1\}$ & $\{1,2\}$ & $\{1,2,3\}$ & $\{1,2,3,4\}$ \\
3       & $\langle 1,2,4,3 \rangle$ &   &    & $\{1,2,4\}$ &      \\
2       & $\langle 1,4,2,3 \rangle$ &   & $\{1,4\}$ &     &      \\
3       & $\langle 1,4,3,2 \rangle$ &   &    & $\{1,3,4\}$ &      \\
2       & $\langle 1,3,4,2 \rangle$ &   & $\{1,3\}$ &     &      \\
1       & $\langle 3,1,4,2 \rangle$ & $\{3\}$ &    &     &      \\
2       & $\langle 3,4,1,2 \rangle$ &   & $\{3,4\}$ &     &      \\
3       & $\langle 3,4,2,1 \rangle$ &   &    & $\{2,3,4\}$ &      \\
2       & $\langle 3,2,4,1 \rangle$ &   & $\{2,3\}$ &     &      \\
1       & $\langle 2,3,4,1 \rangle$ & $\{2\}$ &    &     &      \\
2       & $\langle 2,4,3,1 \rangle$ &   & $\{2,4\}$ &     &      \\
1       & $\langle 4,2,3,1 \rangle$ & $\{4\}$ &    &     &      \\
\hline
\end{tabular}
\caption[Nested set generation for $R_4$] {
Nested set generation for $R_4$ when 
$S_i=\{1,2,\ldots,i\}$ } \label{R4Fig}
\end{figure}
\end{center}

\proclaim{Lemma 2.2}. {
The reverse of $R_m$ is equal to its complement.}

\begin{proof}
We must show that the sum of the $i$-th element of $R_m$
and the $i$-th element of its reverse is always equal to $m$.
This is true for $R_2$.  Assume inductively that it is true
for $R_{m-1}$.  Let $K' = \#(R_{m-1})$. For $i=1,\ldots,K'$,
the sum of the $i$-th element of $R_m$ and the $i$-th element
of its reverse is one more than the sum of the $i$-th element of
$R_{m-1}$ and the $i$-th element of its reverse.  By induction,
this sum is $1+(m-1)=m$.  The middle part of $R_m$ is $[1\ldots,m-1]$,
and this also has the property that the sum of the $i$-th element
of $[1,\ldots,m-1]$ and the $i$-th element of its reverse is $m$.
This completes the induction and the proof.
\end{proof}

One might ask the question: Does the recursive method 
generate all stepping sequences? The answer is ``no''.  For $m=4$,
we exhausted through all sequences in $\{1,2,3\}^{11}$ and found that
exactly 34 of these were stepping sequences.  Theorem~2.1 
produces
$$[3,2,3,2,1,2,3,2,1,2,1],\qquad  [2,3,2,3,1,2,3,1,2,1,2],$$ 
$$[3,2,3,2,1,2,3,1,2,1,2],\qquad [2,3,2,3,1,2,3,2,1,2,1],$$ 
and their reverses, which accounts for 8 of the 34 stepping sequences.
Observe that if two consecutive moves differ by at least two, then one
could do the moves in reversed order, and the only effect would be
that the sets generated on those two moves would be
interchanged.  In particular, if $[i_1,\ldots,i_K]$ is a stepping
sequence and if $|i_j-i_{j+1}| > 1$, then the sequence obtained by
transposing $i_j$ and $i_{j+1}$ is also a stepping sequence.
We will say that these are related by commutation.
Beginning with the eight stepping sequences coming from 
Theorem~2.1, we can produce 10 others that are related
by commutation. The remaining 16 stepping sequences do not seem to arise
from Theorem~2.1.  These are $[2,1,2,3,2,3,1,2,3,2,1]$
and $[2,3,1,2,3,2,1,2,3,1,2]$, together with the stepping sequences
that can be derived from these two using reverse, complement,
and commutation.

\section{The greedy method}

Consider again the nested set problem.  A greedy person wants to alter
the set of largest possible cardinality, subject to the condition that the altered
set has not been seen before.  We show that the sequence $[i_1,i_2,\ldots]$
that is produced in this manner is the stepping sequence $R_m$
that is defined in (\ref{Rm}).
For example, for $m=4$ the reader can check that the greedy method
produces the sequence $[3,2,3,2,1,2,3,2,1,2,1]$, which is the same as $R_4$.

If $S_0,S_1,\ldots,S_m$ are nested sets with $|S_i|=i$, then
for $i=1,\ldots,m-1$ let $S_i^*$ denote the alteration of $S_i$, \ie,
$$S_i^* = S_{i-1} \cup (S_{i+1} \setminus S_i).$$
(Note that $S_{i-1}\subset S_i^* \subset S_{i+1}$, so one can replace
$S_i$ by $S_i^*$ while maintaining the nested set condition.)
The greedy method can formally be stated as follows. In the pseudocode
below, ${\cal S}$ contains
the sets that have been seen before, and $J$ contains the indices $i$
such that $S_i^*$ has not been seen before.

\begin{center}
\begin{figure}[h]
\begin{tabular}{|l|}
\hline
$ G \OpGets{\ } [\ ]$ \\
$ {\cal S}  \OpGets{\ } \{S_1,\ldots,S_m\}$ \\
$ J \OpGets{\ } \{1,2,\dots,m-1\} $ \\
while $J \ne \emptyset$  do $\{$ \\
$\qquad j \OpGets{\ } \max J$ \\
\qquad $ G \OpGets{\ } G \cup [j] $ \\
$\qquad S_j \OpGets{\ } S_j^*$ \\
$\qquad {\cal S} \OpGets{\ } {\cal S} \cup \{S_j\}$ \\
\qquad $ J \OpGets{\ } \{ i \in [1,\ldots,m-1] : S_i^* \not \in {\cal S} \}$ \\
$ \}$ \\
return $G$ \\
\hline
\end{tabular}
\caption{Greedy method for nested sets}
\end{figure}
\end{center}

\proclaim{Theorem 3.1}. {The greedy method produces $G=R_m$.}

\begin{proof}
The theorem is easily seen to hold for $m=2$, because the
greedy sequence is $[1]$ and $R_2=[1]$. Let $m\ge3$.
We will make an inductive hypothesis that the theorem is true for $m-1$, 
and we will show it is true for $m$.  
Recall
$$R_m = (R_{m-1} + 1) \cup [1,\ldots,m-1] \cup R_{m-1}.$$
Suppose initially $S_1=\{a\}$. The greedy method avoids selecting $i=1$ for
as long as possible, since 1 is the smallest possible index. Thus, $S_1$ is
unchanged for the first part of the greedy method, and at this stage all 
elements of ${\cal S}$
contain $a$.     Thus effectively, the greedy algorithm is working with
$S_2,\ldots,S_m$, always selecting the
maximal index such that $S_i^*$ has not been seen before.  By the inductive
hypothesis, $R_{m-1}$, when applied to the nested sets $S_i'=S_i-\{a\}$,
is applying the exact same greedy strategy for determining its next index.
By throwing $a$ back into these sets, we find that $R_{m-1}+1$ is selecting
according to the greedy strategy.   Thus the first $K$ elements of $G$
are $R_{m-1}+1$, where $K=\#(R_{m-1})$.

At this point ${\cal S}$ contains all subsets of $S_m$ that contain $a$.
Let $q_i$ be the unique element in $S_{i+1} \setminus S_i$. Then
$S_i=\{a,q_1,\ldots,q_{i-1}\}$ for $i=0,1,\ldots,m-1$.
Since $S_j^*$ contains $S_{j-1}$, we see that
$a\in S_j^*$ for all $j > 1$, so $S_j^* \in {\cal S}$.
On the other hand, $S_1^* = \{q_1\}$ is not in ${\cal S}$, since $a \not\in S_1^*$.
Thus the greedy method selects~1.
Now $S_2$ contains $a$, so $S_j^*$ contains $a$ for all $j>2$.
On the other hand, $S_2^*=\{q_1,q_2\}$ is new.
Thus, the greedy method selects 2 as the next index.  Continuing in this
way, we see that the greedy algorithm selects $3, 4,\ldots,m-1$.
So the greedy method continues to agree with $R_m$.

At this stage, $S_1=\{q_1\}, S_2=\{q_1,q_2\}, S_3 = \{q_1,q_2,q_3\},
\ldots, S_{m-1} = \{q_1,\ldots,q_{m-1}\}$, and
$$ {\cal S} = \{S \subset S_m : a \in S \} \cup \{S_1,S_2,\ldots,S_{m-1}\}.$$
Note that $S_{m-1}=S_m \setminus \{a\}$.  Consequently,
$$S \not\in {\cal S} \iff \text{$S \subset S_{m-1}$ and $S\ne S_i$
for $i=1,\ldots,m-1$ }$$

The third part of $R_m$ is $R_{m-1}$, applied to $S_0,S_1,\ldots,S_{m-1}$.
By the inductive hypothesis, $R_{m-1}$ follows the
greedy route of selecting the largest index such that $S_i^*$ has not been
seen. This completes the induction and proves the result.
\end{proof}

The greedy person has a humble cousin who always decides to alter
the set of smallest possible cardinality, subject to the
set never having been seen before. 

\proclaim{Lemma 3.2}. {The reverse of the greedy sequence is the humble sequence.}

\begin{proof}  The greedy strategy applied to the complements of $S_i$
results in the humble strategy applied to $S_i$, because the complement
of $S_i$ is largest when $S_i$ is smallest.
Thus, the complement of the greedy sequence $G$ is the humble sequence $H$. 
By Theorem~3.1, $G=R_m$.
By Lemma~2.2,
$H = \text{complement of $G$} = \text{reverse of $G$}$.
\end{proof}

\section{A for-loop method} \label{loopSec}

We showed that $R_m$ can be generated with a recursion or with a greedy
algorithm.
Now we present a third way to generate $R_m$.  If $c$ is an integer,
define $v_2(c)$ to be the largest integer $v$ such that $2^v$ divides $c$ -- 
this is called the valuation of~$c$ at~2. 
Consider the following program.

\begin{center}
\begin{figure}[h]
\begin{tabular}{l}
$G := [\ ]$ \\
for $c=1$ to $2^{m-1}-1$ do $\{$ \\
\quad Let $v = v_2(c)$ (\ie, $ 2^v || c$)  \\
\quad Let $h$ be the Hamming weight of $c$ \\
\quad $d \OpGets{\ } m - v - h$ \\
\quad $G \OpGets{\ } G \cup [d,d+1,\ldots,d+v]$ \\
$\}$  \\
return $G$
\end{tabular}
\caption{For-$c$ loop method to generate $R_m$} \label{forcFig}
\end{figure}
\end{center}

\proclaim{Theorem 4.1}.  {
The sequence $G$ that is returned by the above program is
equal to $R_m$.}

\begin{proof}
We use induction on $m$. Recall the formula
$$R_m = (R_{m-1}+1) \cup [1,2,\ldots,m-1] \cup R_{m-1}.$$
The induction hypothesis shows that $R_{m-1}+1$ is generated with the code:

\begin{center}
\begin{tabular}{l}
$G := [\ ]$ \\
for $c=1$ to $2^{m-2}-1$ do $\{$ \\
\quad Let $v = v_2(c)$  \\
\quad Let $h$ be the Hamming weight of $c$ \\
\quad $d' \OpGets{\ } (m-1) - v - h$ \\
\quad $G \OpGets{\ } G \cup [d'+1,(d'+1)+1,\ldots,(d'+v)+1]$ \\
$\}$  \\
\end{tabular}
\end{center}

\noindent By substituting $d=d'+1$, this is equivalent to
\begin{center}
\begin{tabular}{l}
$G := [\ ]$ \\
for $c=1$ to $2^{m-2}-1$ do $\{$ \\
\quad Let $v = v_2(c)$ \\
\quad Let $h$ be the Hamming weight of $c$ \\
\quad $d \OpGets{\ } m - v - h$ \\
\quad $G \OpGets{\ } G \cup [d,d+1,\ldots,d+v]$ \\
$\}$  \\
\end{tabular}
\end{center}

\noindent Thus, the for loop of the original program with $c$ running from 1 to $2^{m-2}-1$
yields $R_{m-1}+1$.

Next, observe that the for-loop iteration with
$c=2^{m-2}$ sets $v=m-2$, $h=1$, $d=1$, and $G \OpGets{\ } G \cup [1,2,\ldots,m-1]$.
So the for-loop from $c=1$ to $2^{m-2}$ yields $(R_{m-1}+1) \cup [1,2,\ldots,m-1]$.

Finally, the inductive hypothesis implies that the final sequence $R_{m-1}$
can be generated with the code
\begin{center}
\begin{tabular}{l}
for $c'=1$ to $2^{m-2}-1$ do $\{$ \\
\quad Let $v = v_2(c')$ \\
\quad Let $h'$ be the Hamming weight of $c'$ \\
\quad $d \OpGets{\ } (m-1) - v - h'$ \\
\quad $G \OpGets{\ } G \cup [d,d+1,\ldots,d+v]$ \\
$\}$  \\
return $G$
\end{tabular}
\end{center}

\noindent By setting $c=2^{m-2}+c'$ and $h=h'+1$, this can be rewritten as
\begin{center}
\begin{tabular}{l}
for $c=2^{m-2}+1$ to $2^{m-1}-1$ do $\{$ \\
\quad Let $v = v_2(c)$ \\
\quad Let $h$ be the Hamming weight of $c$ \\
\quad $d \OpGets{\ } m - v - h$ \\
\quad $G \OpGets{\ } G \cup [d,d+1,\ldots,d+v]$ \\
$\}$  \\
\end{tabular}
\end{center}

Putting these together, we conclude that $R_m$ can be generated with
the for-loop running from $c=1$ to $2^{m-1}-1$.
This completes the induction and the proof.
\end{proof}

As an example, when $m=4$ the for-$c$ loop generates the following:
$$ \underbrace{[3]}_{c=1} \cup \underbrace{[2\ 3]}_{c=2} \cup \underbrace{[2]}_{c=3}
\cup \underbrace{[1\ 2\ 3]}_{c=4} \cup \underbrace{[2]}_{c=5} \cup
\underbrace{[1\ 2]}_{c=6} \cup \underbrace{[1]}_{c=7}. $$
This is exactly $R_4$.

The advantage of the for-loop method of generating $R_m$, as compared to
the recursive formula (\ref{Rm}), is that it requires less memory.

\section{A fourth way to generate $R_m$}

So far we have given three ways to generate $R_m$: by a recursion,
by the greedy method, and with a for loop.
Here we give a fourth method, also using a for loop.
Consider the following code for $m\ge 2$.

\begin{center}
\begin{figure}[h]
\begin{tabular}{l}
$G \OpGets{\ }[m-1]$ \\
$t \OpGets{\ } m-2$ \\
for $j=1$ to $2^{m-2}-1$ do \\
\qquad let $v = v_2(j)$ \\
\qquad $G \OpGets{\ } G \cup [t,t+1,t+2,\ldots,t+v+1] \cup [t+v]$ \\
\qquad $t \OpGets{+} v-1 $ \\
return $G$ \\
\end{tabular}
\caption{For-$j$ loop method to generate $R_m$} \label{forjFig}
\end{figure}
\end{center}

\proclaim{Theorem 5.1}. {
The above for loop generates $G=R_m$. Also, at the
end of the for loop, $t=0$.}

\begin{proof} We use induction on $m$.  The theorem is easily
seen to be true for $m=2$, because the for loop is empty in that case.
Now let $m\ge 3$, and assume the theorem is true for $m-1$.
We will prove it is true for $m$ also.
The induction hypothesis implies that at the end of the following for loop we will have 
$G = R_{m-1}$ and $t=0$:
\begin{center}
\begin{tabular}{l}
$G \OpGets{\ } [m-2]$ \\
$t \OpGets{\ } m-3$ \\
for $j=1$ to $2^{m-3}-1$ do \\
\qquad let $v = v_2(j)$ \\
\qquad $G \OpGets{\ } G \cup [t,t+1,t+2,\ldots,t+v+1] \cup [t+v]$ \\
\qquad $t \OpGets{+} v-1$ \\
return $G$ \\
\end{tabular}
\end{center}

Since $R_m = (R_{m-1} + 1) \cup [1,\ldots,m] \cup R_{m-1}$, the induction
hypothesis implies that $R_m$ can be generated as follows.
\begin{center}
\begin{tabular}{l|l}
1& $G \OpGets{\ } [m-1]$ \\
2& $t' \OpGets{\ } m-3$ \\
3& for $j=1$ to $2^{m-3}-1$ do \\
4& \qquad let $v = v_2(j)$ \\
5& \qquad $G \OpGets{\ } G \cup [t'+1,t'+2,\ldots,t'+v+2] \cup [t'+v+1]$ \\
6& \qquad $t' \OpGets{+} v-1$ \\
7& $G \OpGets{\ } G \cup [1,2,\ldots,m-1]$ \\
8& $G \OpGets{\ } G \cup [m-2]$ \\
9& $t \OpGets{\ } m-3$ \\
10& for $j=1$ to $2^{m-3}-1$ do \\
11& \qquad let $v = v_2(j)$ \\
12& \qquad $G \OpGets{\ } G \cup [t,t+1,t+2,\ldots,t+v+1] \cup [t+v]$ \\
13& \qquad $t \OpGets{+} v-1 $ \\
14& return $G$ \\
\end{tabular}
\end{center}
Moreover, at line 7 we know that $t'=0$ and at line~14 we know $t=0$.
In lines 1--6, set $t=t'+1$: 
\begin{center}
\begin{tabular}{l|l}
1& $G \OpGets{\ } [m-1]$ \\
2& $t \OpGets{\ } m-2$ \\
3& for $j=1$ to $2^{m-3}-1$ do \\
4& \qquad let $v = v_2(j)$ \\
5& \qquad $G \OpGets{\ } G \cup [t,t+1,\ldots,t+v+1] \cup [t+v]$ \\
6& \qquad $t \OpGets{+} v-1$ \\
\end{tabular}
\end{center}
Then at line~7 we know $t=1$.
Lines 7, 8, and 9 are equivalent to
\begin{center}
\begin{tabular}{l}
$j \OpGets{\ }2^{m-3}$ \\
Let $v = v_2(j)$  (\ie, $v=m-3$) \\
$G \OpGets{\ } G \cup [t,t+1,\ldots,t+v+1] \cup [t+v]$ \\
$t \OpGets{+} v-1$ \\
\end{tabular}
\end{center}
We recognize this as the iteration of the for loop with $j=2^{m-3}$.
Finally, in the last for loop we can change $j$ to $j+2^{m-3}$ without
affecting the value for $v$. Equivalently, line~10 can be changed
to ``for $j=2^{m-3}+1$ to $2^{m-2}-1$ do''.
This shows that $R_m$ can be produced
by the above for-loop with $j$ running from 1 to $2^{m-2}-1$.
This completes the induction and the proof.
\end{proof}

We remark that the two for-loop methods for generating $R_m$ that
are given in Figures~\ref{forcFig} and~\ref{forjFig} are closely related,
as shown in the next lemma.

\proclaim{Lemma 5.2}.  {The update for $G$ in the for loop of Figure~\ref{forjFig}
is equivalent to the updates when $c=2j$ and $c=2j+1$ in the for loop
of Figure~\ref{forcFig}. }

For example, when $m=4$ the for-$c$ loop of Theorem~4.1
and the for-$j$ loop of Theorem~5.1
generate $R_4$ as follows:
$$ \underbrace{[3]}_{c=1} \cup 
\underbrace{\underbrace{[2\ 3]}_{c=2} \cup \underbrace{[2]}_{c=3}}_{j=1}
\cup 
\underbrace{\underbrace{[1\ 2\ 3]}_{c=4} \cup \underbrace{[2]}_{c=5}}_{j=2}
\cup
\underbrace{\underbrace{[1\ 2]}_{c=6} \cup \underbrace{[1]}_{c=7}}_{j=3}. $$

\begin{proof}
Let $v_j$ denote the valuation of~$v$ at~2, (\ie, $2^{v_j}$ maximally divides~$j$.
Let $t_j$ denote the value for $t$ at the point in the for-$j$ loop 
when $G$ is being updated.  Let $d_c$ denote the value for $d$ in
the for-$c$ loop, namely $d_c=m-v_c-\HW(c)$.  The for-$j$ loop does the update 
\begin{equation}
G \OpGets{\ } G \cup [t_j,t_j+1,\ldots,t_j+v_j+1] \cup [t_j+v_j] 
\label{forjUpdate}
\end{equation}
and the for-$c$ loop does the update
$$G \OpGets{\ } G \cup [d_c,d_c+1,\ldots,d_c+v_c].$$
We claim that when $c=2j$ then 
$$[d_c,\ldots,d_c+v_c] = [t_j,\ldots,t_j+v_j+1]$$
and when $c=2j+1$ then
$$[d_c,\ldots,d_c+v_c] = [t_j+v_j].$$
The claim will imply that the update (\ref{forjUpdate}) in the for-$j$
loop is equivalent to the two updates with $c=2j$ and $c=2j+1$ in the
for-$c$ loop.
To see the claim, it suffices to show that $d_{2j}=t_j$, $v_{2j}=v_j+1$,
$d_{2j+1}=t_j+v_j$, and $v_{2j+1}=0$. The statements $v_{2j}=v_j+1$
and $v_{2j+1}=0$ are obvious, so we just need to prove
$$\text{$d_{2j}=t_j$ and $d_{2j+1}-d_{2j}=v_j$.}$$
Looking at the for-$j$ loop program,
we see that $t$ starts at $m-2$ and then $v_j-1$ is added to $t$
at the end of each for-loop iteration.  Thus,
$$t_j+v_j-1= m-2+\sum_{0<i\le j} (v_i-1).$$
Lemma~5.3 below shows that the right-hand side is 
$m-2-\HW(j)$, and so
$$t_j= m-2-\HW(j)-v_j+1=m-1-\HW(j)-v_j.$$
On the other hand,
$d_{2j}=m-v_{2j}-\HW(2j)=m-v_j-1-\HW(j)$, and so we see that $d_{2j}=t_j$
as claimed. Next, 
\begin{eqnarray*}
d_{2j+1}-d_{2j} &=& (m-v_{2j+1}-\HW(2j+1))-(m-v_{2j}-\HW(2j)) \\
&=& v_{2j}+\HW(2j)-\HW(2j+1) = v_{2j}-1=v_j.
\end{eqnarray*}
This proves the claim and establishes the lemma.
\end{proof}

\proclaim{Lemma 5.3}. {
If $i>0$, let $v_i$ denote the valuation of~$i$ at~2 (\ie, $2^{v_i}$ maximally divides~$i$), 
and let $\HW(i)$ denote the Hamming weight of~$i$. Then
$$ \sum_{i=1}^j (1-v_i) = \HW(j).$$ }

\begin{proof}  This is trivially true when $j=1$.  Now let $j\ge 1$.
Assuming the lemma is true for $j$, we will show it is true for $j+1$.
By the inductive hypothesis,
$$\sum_{i=1}^{j+1} (1-v_i) = \HW(j) + 1 - v_{j+1}$$
and we must show that this is equal to $\HW(j+1)$.  That is, we  must show
that
\begin{equation}
\HW(j+1)-\HW(j) = 1 - v_{j+1}. \label{hwEqn}
\end{equation}
First, if $j$ is even then $j+1$ is obtained by changing the low bit of
$j$ from a 0 to a 1.  Thus, both sides of the above equality are 1.
Next, if $j$ is odd, then let $k$ be the number of consecutive 1's in
the low bits of $j$. Then $j+1$ is obtained from $j$ by complementing the lowest
$k+1$ bits.  That is, the bits $011\ldots 1$ in $j$ become $100\ldots0$
in $j+1$.  This shows that $v_{j+1}=k$ and $\HW(j)-\HW(j+1)=k-1$,
and (\ref{hwEqn}) easily follows. This completes the induction and the proof.
\end{proof}

Mark Jacobson observed that $\sum_{i=1}^j v_i$ is the valuation
of~$j!$ at~2.  Thus, Lemma~5.3 is equivalent to
the statement that the valuation of~$j!$ at~2 is
$j-\HW(j)$.  This latter fact is proved in \cite[Section 4.4]{GKP},
using a different proof from the above.

\section{Relationship with other combinatorial Gray codes}

In this section we show that 
the stepping sequences $R_m$ give rise to a new family of cyclic Gray codes
through binary $m$-bit integers.  Also they give rise to new 
Gray codes through $k$-element subsets of an $m$-element set.

We first explain the connection with Gray codes on binary $m$-bit integers.
If $S_m=\{0,1,\ldots,m-1\}$, then subsets of $S_m$ are
in bijection with the integers in the range $0 \le i < 2^m$ by
the map $S \mapsto \sum_{i \in S} 2^i$.  The cardinality of $S$
is equal to the Hamming weight of its corresponding integer.
When we Gray-code through nested subsets of $S_m$, we exhaust through all of
its subsets, and consequently we exhaust through the $m$-bit integers.
For example, the stepping sequence $R_4$ exhausts through subsets of
$\{0,1,2,3\}$ in the following order (see Figure~\ref{R4Fig}):
$$\emptyset, \{0\}, \{0,1\}, \{0,1,2\}, \{0,1,2,3\}, \{0,1,3\},
\{0,3\}, \{0,2,3\}, $$
$$\{ 0,2\}, \{2\}, \{2,3\}, \{1,2,3\}, \{1,2\}, \{1\}, \{1,3\}, \{3\}.$$
This corresponds the sequence of integers (written in binary):
$$0000, 0001, 0011, 0111, 1111, 1011, 1001, 1101, 0101, 0100, 1100, 1110,
0110, 0010, 1010, 1000.$$
In particular, this is a Gray-code ordering of the 4-bit integers,
because at each step, only one bit is altered. 
A Gray code is {\it cyclic} if the last integer in the
sequence differs by only one bit from the first integer in the sequence.
The above example is a cyclic Gray code, because the first element is 0000 
and the last element has Hamming weight~1.

It is not true that
all stepping sequences give rise to Gray-code orderings of integers.
For example, the stepping sequence
$G=[3,2,3,2,1,2,3,1,2,1,2]$ produces an integer of Hamming weight~3
next to an integer of Hamming weight~1.  

We will say that a sequence of integers is {\it contiguous} if 
consecutive elements differ by $\pm1$.  We will say that a stepping
sequence for $m$ is {\it strongly contiguous} if it is contiguous,
it begins with $m-1$, and it ends with 1.

\proclaim{Lemma 6.1}. {A stepping sequence gives rise to a cyclic Gray code
if and only if it is strongly contiguous.}

\begin{proof}  
Let $G=[g_0,\ldots,g_K]$ be a stepping sequence.
The initial sets $S_0,S_1,\ldots,S_m$ correspond to the integers
$$00\cdots000, 00\cdots001, 00\cdots011, \ldots, 01\cdots111, 11\cdots111.$$ 
The next integer produced by $G$
has Hamming weight equal to $g_0$.  This integer is one-off from
$1\cdots1$ if and only if $g_0=m-1$.
For $0\le i < K$, the integer produced by the $i$-th and $(i+1)$-th
elements of $G$ have Hamming weights $g_i$ and $g_{i+1}$. 
Because the nesting
property of the sets is always maintained, we know that the one-bits of 
the integer
with lower Hamming weight is a subset of the one-bits of the integer of higher
Hamming weight.  Thus, these two integers differ by a single bit if and
only if their Hamming weights differ by one, \ie, $|g_i-g_{i+1}|=1$.
Finally, the cyclic condition
mandates that the last element produced be adjacent to $0\cdots000$, 
{\it i.e.}, the last element must have Hamming weight~1. This is equivalent
to $g_K=1$.
\end{proof}

\proclaim{Theorem 6.2}. {
The stepping sequence $R_m$ is strongly contiguous for each $m\ge2$, 
and so it gives rise to a cyclic Gray-code ordering on the $m$-bit integers.
More generally, if $A$ and $B$ are strongly contiguous stepping sequences for $m-1$, then
$C = (A+1) \cup [1,2,\ldots,m-1] \cup B$ is a strongly contiguous stepping sequence
for $m$. }

\begin{proof}  We first prove the second statement. Note that $C$ is a stepping
sequence for $m$ by Theorem~2.1, and we just need to show that it is
strongly contiguous.  Since $A$ and $B$ are contiguous, so are the three pieces
that make up $C$.  Also, since $A$ ends with a~1, we know $A+1$ ends with a 2, which
is one-off from the first term of $[1,2,\ldots,m-1]$. Likewise, $B$ begins with $m-2$,
which is one-off from the last term of $[1,2,\ldots,m-1]$. This shows that $C$ is contiguous.
The first term of $A$ is $m-2$, so the first term of $A+1$ is $m-1$. Also, the last term of $B$ is~1.
Thus, the first term of $C$ is $m-1$ and the last term of $C$ is~1. This shows that $C$
is strongly contiguous.

Taking $A=B=R_{m-1}$, we see that if $R_{m-1}$ is strongly contiguous then so is $R_m$.
Since $R_2=[1]$ is strongly contiguous for $m=2$, it follows by induction that $R_m$ is strongly 
contiguous for all~$m\ge 2$.
\end{proof}

For $m=4$ and $m=5$, an exhaustive search shows that $R_m$ is the only strongly contiguous 
stepping sequence.  However, for $m=6$ there are exactly two strongly
contiguous stepping sequences: $R_6$ and $A_6$, where
\begin{eqnarray*}
A_6 &=& [ 5, 4, 5, 4, 3, 2, 3, 4, 5, 4, 3, 2, 3, 4, 3, 2, 3, 2, 1, 2, 3, 
4, 5, 4, 3, 4, 3, 2,  \\
&& 3, 4, 3, 2, 3, 2, 1, 2, 3, 4, 5, 4, 3, 4, 3, 2, 3, 4, 3, 2, 1, 2, 
    3, 4, 3, 2, 1, 2, 1 ]. 
\end{eqnarray*}
An interesting observation is that $A_6$ and $R_6$ have a long
subsequence in common:
\begin{eqnarray*}
A_6 &=& [ 5, 4, 5, 4, 3, 2, 3, 4, 5, \underline{4, 3, 2, 3, 4, 3, 2, 3, 2, 1, 2, 3, 
4, 5, 4, 3, 4, 3, 2,}  \\
&& \underline{3, 4, 3, 2, 3, 2, 1, 2, 3, 4}, 5, 4, 3, 4, 3, 2, 3, 4, 3, 2, 1, 2, 
    3, 4, 3, 2, 1, 2, 1 ] \\ 
R_6 &=& [ 5, 4, 5, 4, 3, 4, 5, 4, 3, 4, 3, 2, 3, 4, 5, 4, 3, \underline{4, 3, 2, 3, 
4, 3, 2, 3, 2, 1, 2,}  \\
&&\underline{3, 4, 5, 4, 3, 4, 3, 2, 3, 4, 3, 2, 3, 2, 1, 2, 3, 4}, 3, 2, 3, 2, 1, 2, 3,
2, 1, 2, 1 ]
\end{eqnarray*}
Notice that the reverse of $A_6$ is equal to its complement.  However, it is not
true that all strongly contiguous stepping sequences have that property.
The first counterexamples occur when $m=7$:
$(A_6+1) \cup [1,\ldots,6]
\cup R_6$ and $(R_6+1) \cup [1\ldots,6] \cup A_6$.

Theorem~6.2 implies that if there are $N$ strongly contiguous stepping
sequences for $m-1$, then there are at least $N^2$ strongly contiguous stepping sequences for~$m$, 
because when constructing $C=(A+1)\cup[1,\ldots,m-1]\cup B$ we have 
$N$ choices for~$A$ and $N$ choices for~$B$.
Thus, if $c_m$ denotes the number of strongly contiguous stepping sequences 
for $m$, then
$c_m \ge c_{m-1}^2$.  By exhaustion we found $c_m=1$ for $m\le 5$, and $c_6=2$. Thus,
$c_7 \ge 4$, $c_8 \ge 16$, and in general $c_m \ge 2^{2^{m-6}}$ for $m\ge 6$.
For $m \ge 7$, it is an open question whether $c_m =2^{2^{m-6}}$.

By Lemma~6.1, each strongly contiguous stepping sequence gives rise
to a Gray code on $m$-bit integers. Thus, we have an explicit recipe for
constructing $c_m$ new Gray codes on $m$-bit integers.  As shown above,
the number of these is at least $2^{2^{m-6}}$.

One might wonder how many stepping sequences there are that are
contiguous but not strongly contiguous.
By exhaustive computer search we found that for $m\le 6$, every contiguous stepping sequence
has the property that either the sequence or its reverse is strongly contiguous. It is
an open question whether this is true for larger~$m$. The only partial result that we have in this
direction is given in Lemma~6.4 below.

\proclaim{Lemma 6.3}. {If $m$ is even, then 
a contiguous stepping sequence must begin and end with an odd integer.}

\begin{proof}
Let $m$ be even. Note that $R_m$ begins with the odd integer $m-1$ and 
ends with the odd integer~1. 
Since every other step produces a set of even
cardinality and every other step produces a set of odd cardinality, and since the total number
of steps is the odd number $2^m-m-1$, we see that if one begins on a set of even cardinality,
then more even-order sets are produced than odd-order sets.  However, $R_m$ is known to be
a stepping sequence that begins on an odd integer.  Then the odd-order sets (excluding the
initial ones) must outnumber the even-order sets. It follows that any contiguous stepping sequence
must begin with an odd integer.
\end{proof}

\proclaim{Lemma 6.4}.  {Let $S$ be a contiguous stepping sequence. Then $S$ or 
its reverse has the property that its first element is congruent to
$m-1 \pmod 2$ and the last element is congruent to $1 \pmod 2$.}

\begin{proof}
When $m$ is even, the result follows immediately from Lemma~6.3.
When $m$ is odd, then $2^m-m-1$ is even, and since the parities of the
elements of $S$ alternate, we conclude that the first and last elements
have opposite parities.  By reversing the 
stepping sequence if necessary, we can assume that the first element is
even and the last is odd.  The result follows.
\end{proof}

We showed that Gray codes through nested sets
give rise to a Gray code through $m$-bit integers, provided the
associated stepping sequence is contiguous.
However, it is {\it not} true that a Gray code through the $m$-bit integers 
gives rise to
a solution to the nested set problem.  For example, consider
the binary reflective Gray code sequence that is defined recursively by
$G_1=[0,1]$ and $G_m = ([0] \x G_{m-1} ) \cup ([1] \x G_{m-1}')$, where
$G_{m-1}'$ denotes the reverse of $G_{m-1}$.  For $m=1,2,3$
these are
$$G_1=[0,1],\qquad G_2 = [00,01,11,10],\qquad
G_3=[000,001,011,010,110,111,101,100].$$
The first few corresponding sets for $G_3$ are
\begin{center}
\begin{tabular}{c c c c}
$\emptyset$ & & &  \\
& \{0\}  &&  \\
&& \{0,1\}  &  \\
& \{1\}  &&  \\
&& \{1,2\}  &  \\
&& &\{0,1,2\}    \\
&& \{0,2\}  &  \\
\end{tabular}
\end{center}
and at this point the sets $\emptyset, \{1\}, \{0,2\}, \{0,1,2\}$ are not
nested.

We close this section by pointing out a connection with the problem
of Gray-coding through $k$-element subsets of an $m$-element set.
A Gray code through $k$-element subsets is defined as an ordering of
the $k$-element subsets such that two consecutive subsets differ by
just one element, \ie, they have exactly $k-1$ elements in common.
This ordering is {\it cyclic} if the first and last $k$-element subsets
also differ by just one element.
A Gray code through nested sets results in an ordered sequence of all subsets
of $S_m$, that contains each subset exactly once.  We claim that by restricting 
to the $k$-element subsets, we obtain a Gray code through $k$-element subsets.
To see this, recall that a new $k$-element set is produced when 
$S_k$ is replaced by a new set $S_k^*$.  To obtain $S_k^*$ from $S_k$,
the unique element of $S_k\setminus S_{k-1}$ is
replaced by the unique element of $S_{k+1}\setminus S_k$.  Thus, the
two sets differ by a single element, and this demonstrates
that a Gray code through nested sets induces a Gray code through $k$-element
subsets.

For example, the stepping sequence $R_6$ induces
the following order on 2-element subsets of $\{0,1,2,3,4,5\}$: 
$$\{0,1\}, \{0,5\}, \{0,4\}, \{0,3\}, \{0,2\}, \{2,3\}, \{2,5\}, \{2,4\}, $$
$$\{1,2\}, \{1,4\}, \{1,3\}, \{1,5\}, \{3,5\}, \{4,5\}, \{3,4\}.$$
Note however that the first and last sets are not close (they are disjoint), so this
combinatorial Gray code through order-2 subsets is not cyclic.
For methods of Gray-coding through $k$-element subsets of an order-$m$ set see \cite{Chase,Wilf,BW}.

\end{document}